\def\pdual{\mathop{\hbox{\rm $p$-dual}}\nolimits}

\documentclass[12pt]{amsart}
\usepackage{amssymb}

\def\above#1#2#3#4{\genfrac{}{}{0pt}{}{\lower#1\hbox{\smash{$\scriptstyle #3$}}}{\raise#2\hbox{\smash{$\scriptstyle #4$}}}}

\newdimen{\twowidth}
\newdimen{\halftwo}
\def\slashtwo{2\llap{\hbox to 0pt{\hss$|$\hss}\kern\halftwo}}
\newdimen{\threewidth}
\newdimen{\halfthree}
\def\slashthree{3\llap{\hbox to 0pt{\hss$|$\hss}\kern\halfthree}}
\newdimen{\inftywidth}
\newdimen{\halfinfty}
\def\slashinfty{\infty\llap{\hbox to 0pt{\hss$|$\hss}\kern\halfinfty}}

\def\BreakWithGlue#1#2#3{\nolinebreak\hskip #1\hbox{}\penalty #2\hbox{}\nolinebreak\hskip #3}
\def\EasyButWeakLineBreak{\BreakWithGlue{0pt plus 1 fill}{500}{0pt plus 1 fill}}
\def\HardButStrongLineBreak{\BreakWithGlue{0pt plus 1 filll}{600}{0pt plus -1 filll}}
\newcommand\spacer{\ }
\def\instructions#1{$\langle #1\rangle$}
\def\pdualwitharg#1{\mathop{\hbox{$#1$-dual}}\nolimits}
\def\pfillwitharg#1{\mathop{\hbox{$#1$-fill}}\nolimits}
\def\main{\mathop{\rm main}\nolimits}
\def\iso{\cong}

\def\NumLattices{8595}

\begin{document}

\title[Unabridged table of reflective lattices]{Unabridged table of reflective Lorentzian lattices of rank 3}
\author{Daniel Allcock}
\address{Department of Mathematics\\U.T. Austin}
\email{allcock@math.utexas.edu}
\urladdr{http://www.math.utexas.edu/\textasciitilde allcock}
\thanks{Partly supported by NSF grant DMS-0600112.}
\date{November 19, 2010}
%\date{November 14, 2010}

\begin{abstract}
This is the unabriged table of all \NumLattices{} rank three
reflective Lorentzian lattices, intended as a supplement to the
author's paper classifying them.  The abridged table in that paper is
complete too, but less explicit, so the main purpose of this document
is archival.  The TeX sourcecode is simultaneously a Perl script,
which when run prints out all the lattices in computer-readable
format.
\end{abstract}

\maketitle

\noindent
See the author's paper \cite{Allcock-rk3} for the meaning of this
table, as well as the proof and more detailed references.  The
following comments describe features of the table present here but not
in the paper.

{\it Names and constructions of Lattices:\/} For a given Weyl group
$W_n$, we index the lattices having that Weyl group as in \cite{Allcock-rk3},
followed by those that are printed implicitly there, followed by those
that are got by $p$-duality operations.  For the lattices printed here but not there,
the method of construction is stated.  When $L$ is printed explicitly
or implicitly in the paper, and $M$ is got from it by duality, we have
printed all chains of $p$-dualities from $L$ to $M$.  (Except that we
omit chains involving a prime $p$ if $L\iso\pdual(L)$ or equivalently
$M\iso\pdual(M)$.)  This means that some lattices have more than one
construction; no one of those constructions is
distinguished, but the set of all of them is natural.

{\it Reference to Nikulin's lattices:\/}  When one of our lattices appears on
Nikulin's tables \cite{Nikulin-rk3}, we indicate where it does so.

{\it Gram matrices:\/} We computed the Gram matrices and typeset them,
to see that they really do coincide with Nikulin's Gram
matrices. Unfortunately this made the TeX file too large for the arxiv
to accept.  If you want them they you can download a version of these
tables containing the Gram matrices from {\it
http://www.math.utexas.edu/$\sim$allcock} (the author's web page).  Or you can
extract the computer-readable data and compute them
yourself.

{\it Extracting computer-readable data:\/}  By arcane methods we have
arranged for the TeX source to be simultaneously a Perl script that
prints out all \NumLattices{} lattices (or just some of
them) in a computer-readable format, namely the one read by the
PARI/GP software.  If the file is saved as {\tt
file.tex} then simply enter {\tt perl file.tex} at the unix command
line and follow the instructions that will appear.

\vfill\eject
\begin{center}
\bf Unabridged Table of Rank 3 Reflective Lorentzian Lattices
\end{center}
\medskip
\begingroup
\parindent=0pt
\parskip=0pt\small
\raggedbottom

% for shrinking to arxiv-friendly size
\def\A{\leavevmode\llap{\spacer}}
\def\B{{$\left[\!\llap{\phantom{\begingroup \smaller\smaller\smaller\begin{tabular}{@{}c@{}}0\\0\\0\end{tabular}\endgroup}}\right.$}}
\def\C{{$\left.\llap{\phantom{\begingroup \smaller\smaller\smaller\begin{tabular}{@{}c@{}}0\\0\\0\end{tabular}\endgroup}}\!\right]$}}
\def\D{\begingroup \smaller\smaller\smaller\begin{tabular}{@{}c@{}}}
\def\E{\end{tabular}\endgroup\kern3pt}
\def\F{\nopagebreak\par}
\def\G#1{\spacer\instructions{#1}}
\def\H{\EasyButWeakLineBreak}
\def\I{\nopagebreak\par\leavevmode}
\def\J{\end{tabular}\endgroup}
\def\K{\HardButStrongLineBreak\kern3pt}
\def\L{\hbox{}\par\smallskip}
\def\M#1#2{\above{1pt}{1pt}{#1}{#2}}
\def\N{\J\K\D}
\end{document}